\documentclass{amsart}[12pt]
\usepackage{amsfonts,amssymb,amsmath,amsthm}
\usepackage{enumerate}
\usepackage{amssymb,latexsym}
\usepackage{mathrsfs}

\newtheorem{rem}{Remark}
\usepackage{fullpage}
\title{A stationary process associated with the Dirichlet distribution arising from the complex projective space}
\author[N. Demni]{N. Demni}
\address{Institut de Recherche en Math\'ematiques de Rennes\\ Universit\'e Rennes 1\\
France}
\email{nizar.demni@univ-rennes1.fr}
\begin{document}
\maketitle
\begin{abstract}
Let $(U_t)_{t \geq 0}$ be a Brownian motion valued in the complex projective space $\mathbb{C}P^{N-1}$. Using unitary spherical harmonics of homogeneous degree zero, we derive the densities of $|U_t^{1}|^2$ and 
of $(|U_t^{1}|^2, |U_t^2|^2)$, and express them through Jacobi polynomials in the simplices of $\mathbb{R}$ and $\mathbb{R}^2$ respectively. More generally, the distribution of $(|U_t^{1}|^2, \dots, |U_t^k|^2), 2 \leq k \leq N-1$ may be derived using the decomposition of the unitary spherical harmonics under the action of the unitary group $\mathcal{U}(N-k+1)$ yet computations become tedious. We also revisit the approach initiated in \cite{Nec-Pel} and based on a partial differential equation (hereafter pde) satisfied by the Laplace transform of the density. When $k=1$, we invert the Laplace transform and retrieve the expression derived using spherical harmonics. For general $1 \leq k \leq N-2$, the integrations by parts performed on the pde lead to a heat equation in the simplex of $\mathbb{R}^k$.
\end{abstract}
\section{Motivation}
The complex unit sphere 
\begin{equation*}
S^{2N-1} :=  \{(z_1,\dots, z_N), |z_1|^2 + \dots + |z_N|^2 = 1\},  \quad N \geq 1,
\end{equation*} 
is a compact manifold without boundary therefore carries a Brownian motion $(U_t)_{t \geq 0}$ defined by means of its Laplace-Beltrami operator. This process is stationary and the random variable $U_t$ converges weakly as $t \rightarrow \infty$ to a uniformly-distributed random vector $U_{\infty}$. For the latter, it is already known that 
\begin{equation*}
(|U_{\infty}^{1}|^2, \dots, |U_{\infty}^{k}|^2), \quad 1 \leq k \leq N-1,
\end{equation*}
follows the Dirichlet distribution (\cite{Hia-Pet})
\begin{equation}\label{Diri}
s_k(u) \prod_{i=1}^kdu_i \triangleq (1-u_1-u_2-\dots-u_k)^{N-k-1}{\bf 1}_{\Sigma_k}(u)\prod_{i=1}^kdu_i
\end{equation}
where $\Sigma_k = \{u_i > 0, 1 \leq i \leq k, u_1+\dots + u_k < 1\}$ is the standard simplex. Motivated by quantum information theory, the investigations of the distribution of 
\begin{equation*}
U_t^{(k)} \triangleq (|U_{t}^{1}|^2, \dots, |U_{t}^{k}|^2), \quad 1 \leq k \leq N
\end{equation*}
started in \cite{Nec-Pel} yet have not been completed. There, a linear pde for the Laplace transform of this distribution was obtained and partially solved only when $k=1$. Recall that for the Brownian motion on the Euclidian sphere $S^{N-1}$, the density of a single coordinate is given by a series involving products of ultraspherical polynomials of index $(N-2)/2$ (\cite{Kar-McG}). The main ingredients leading to this series are the expansion of the heat kernel on $S^{N-1}$ in the basis of $O(N)$-spherical harmonics and on Gegenbauer addition Theorem (\cite{Wat}, p.369). In the complex setting, it is very likely known that $(|U_t^{1}|^2)_{t \geq 0}$ is a real Jacobi process (see \cite{Bak} and references therein). Nonetheless, one wonders how does the proof written in \cite{Kar-McG} carry to the Brownian motion on $S^{2N-1}$ and how does it extend in order to derive the density of $U_t^{(k)}$. In the first part of this paper, we answer these questions by considering the heat kernel on the complex projective space $\mathbb{C}P^{N-1} = S^1/S^{2N-1}$ rather than $S^{2N-1}$. This is by no means a loss of generality since we are interested in the joint distribution of the moduli of $k$ coordinates of $U_t$. Besides, the space of continuous functions on $\mathbb{C}P^{N-1}$ decomposes as the direct sum of subspaces of $\mathcal{U}(N)$-spherical harmonics that are homogeneous of degree zero, while the decomposition of continuous functions on $S^{2N-1}$ involves all spherical harmonics (\cite{Gri}). Accordingly, the heat kernel on $\mathbb{C}P^{N-1}$ is expressed as a series of normalized Jacobi polynomials $(P_n^{N-2,0}/P_n^{N-2,0}(1))_{n \geq 0}$ which for each $n$, gives the $n$-th reproducing kernel on $\mathbb{C}P^{N-1}$ (\cite{Koo}). Hence, the integration over the sphere $S^{2N-3}$ together with an application of Koornwinder's addition Theorem (\cite{Koo}) lead to the density of $|U_t^{1}|^2$: 
\begin{align}\label{JacDen}
f_t(c,u) &\triangleq \left[\sum_{n = 0}^{\infty} e^{-n(n+N-1)t}\frac{P_n^{N-2,0}(2c-1)P_n^{N-2,0}(2u-1)}{||P_n^{N-2,0}||_2^2}\right] s_1(u)
\end{align}
where we set $c \triangleq |U_0^{1}|^2 \in [0,1]$ and $||P_n||_2^2$ is the squared $L^2$-norm of $u \mapsto P_n^{N-2,0}(2u-1)$ with respect to $s_1(u)du$. 

Up to an additional ingredient, the derivation of the density of $U_t^{(2)}$ is quite similar. Loosely speaking, we would like to integrate the heat kernel over the sphere $S^{2N-5}$ (we assume $N$ large enough) and as such, we need to decompose degree zero homogenous spherical harmonics in $S^{2N-1}$ under the action of the unitary group $\mathcal{U}(N-1)$. This decomposition is stated in \cite{Koo1}, Theorem 5.1, and the $n$-th reproducing kernel in turn decomposes as a wighted sum of reproducing kernels on $S^{2N-3}$. Consequently, Koornwinder's addition Theorem again leads to the sought density which may be expressed through Jacobi polynomials in the simplex $\Sigma_2$ (\cite{Dun-Xu}, Proposition 2.3.8 p.47). More precisely, if we denote these polynomials by $(Q_{j,n-j}^{(N)})_{n \geq 0 , 0 \leq j \leq n}$ then the density reads 
\begin{equation}\label{JacDen1}
\left[\sum_{n=0}^{\infty} e^{-n(n+N-1)t} \sum_{j=0}^n\frac{Q_{n-j,j}^{(N)}(c_1,c_2)Q_{n-j,j}^{(N)}(u_1,u_2)}{||Q_{n-j,j}^{(N)}||_2^2}\right]s_2(u),
\end{equation}
where we set $(c_1,c_2) \triangleq (|U_0^1|^2, |U_0^2|^2) \in \Sigma_2$ and $||Q_{n-j,j}^{(N)}||_2^2$ is the squared $L^2$-norm of $Q_{n-j,j}^{(N)}$ with respect to $s_2(u)du_1du_2$.

More generally, the derivation of the density of $U_t^{(k)}, 2 \leq k \leq N-1$ relies on the decomposition of the spherical harmonics under the action of $\mathcal{U}(N-k+1)$ and is expressed through orthonormal Jacobi polynomials in $\Sigma_k$ as 
\begin{equation*}
\sum_{n \geq 0}e^{-n(n+N-1)t} \sum_{\tau \in \mathbb{N}^k, |\tau| = n} Q_{\tau}^{(N)}(c_1,c_2,\dots, c_k)Q_{\tau}^{(N)}(u)s_k(u)
\end{equation*}
where $(c_1, \dots, c_k) \triangleq (|U_0^1|^2, \dots, |U_0^k|^2)$. Yet computations become tedious and we are not willing to exhibit them here. Rather, we shall revisit and complete the investigations started in \cite{Nec-Pel}. Actually, an expression for the Laplace transform of the density of $|U_t^1|^2$ was obtained there and involves the following sequence $(a_n = a_n(c,N))_{n \geq 0}$ of real numbers determined recursively by (\cite{Nec-Pel}, eq. 4.23)
\begin{equation}\label{S1}
\sum_{n=0}^p a_n \binom{p}{n} \frac{1}{(N+2n)_{p-n}} = \frac{c^p}{p!}, \quad p \geq 0, \quad a_0 = 1.
\end{equation}
In particular, the following was proved (\cite{Nec-Pel} eq. 4.24. and eq. 4.25):
\begin{equation*}
a_n(0,N)  = \frac{(-1)^n}{(N+n-1)_n}, \quad a_n(1,N) = \frac{(N-1)_n}{n!(N+n-1)_n},
\end{equation*}
which we can rewrite as 
\begin{equation*}
\frac{1}{(N+n-1)_n}P_n^{N-2,0}(-1), \quad \frac{1}{(N+n-1)_n} P_n^{N-2,0}(1)
\end{equation*}
respectively. Using a Neumann series for Bessel functions (\cite{Wat}), we shall prove that for all $c \in [0,1]$
\begin{equation}\label{sequence}
a_n = a_n(c,N) = \frac{1}{(N+n-1)_n} P_n^{N-2,0}(2c-1). 
\end{equation}
Having these coefficients in hands, we can then invert the Laplace transform and retrieve \eqref{JacDen}. At this level, we point out that the pde satisfied by the Laplace transform of the density of $|U_t^1|^2$ leads after integrations by parts to the heat equation associated with the Jacobi operator
\begin{equation}\label{JacOp}
u(1-u)\partial_u^2 + [1 - Nu]\partial_u. 
\end{equation}
More generally, the pde satisfied by the Laplace transform of the joint distribution of $U_t^{(k)}, 1 \leq k \leq N-2$ gives rise to the heat equation on the standard simplex associated with the generalized Jacobi operator (\cite{Akt-Xu}, see also \cite{Dun-Xu} p.46 but consult the list of errata available on the webpage of Y. Xu):  
\begin{equation}\label{HE2}
\sum_{i=1}^k[1 - N u_i]\partial_i+ \sum_{i=1}^k(u_i-u_i^2)\partial_{ii} - \sum_{i \neq j}u_iu_j\partial_{ij}.
\end{equation}
This is an elliptic operator admitting different orthogonal basis of eigen-polynomials corresponding to the sequence of eigenvalues $\{-n(N+n-1), n \geq 0\}$. Among them figure the Jacobi polynomials in the simplex which agrees with our previous computations. 

The paper is organized as follows. The two following sections are concerned with the derivations of the densities of $U_t^{(k)}, k\in \{1,2\}$. In section $4$, we solve the system \eqref{S1} and invert the Laplace transform of the density of $|U_t^1|^2$. In section $5$, we perform integrations by parts on the pde satisfied by the Laplace transform of the density of $U_t^{(k)}$, omitting for a while the boundary terms. In the last section, we write down the latters and show that all of them vanish unless $k = N-1$.
\section{The distribution of $|U_t^1|^2$} 
Let $m,n$ be non negative integers and recall from \cite{Koo} that $(m,n)$-complex spherical harmonics are the restriction $S^{2N-1}$ of harmonic polynomials in the variables 
\begin{equation*}
(z_1, z_2, \dots, z_N, \overline{z_1}, \overline{z_2}, \dots, \overline{z_N})
\end{equation*}
which are $m$-homogenous in the variables $(z_i)_{i = 1}^N$ and $n$-homogeneous in the variables $(\overline{z_i})_{i = 1}^N$. Taking $m=n$, we obtain the $(n,n)$-complex spherical harmonics that are homogenous of degree zero with respect to the action of $S^1$. Their restrictions to $\mathbb{C}P^{N-1}$ form a dense algebra in the space of continuous functions on $\mathbb{C}P^{N-1}$ endowed with to the uniform norm (see \cite{Gri}, p.189). Moreover, the spectrum of the Laplace-Beltrami operator on $\mathbb{C}P^{N-1}$ is given by the sequence $\{-n(n+N-1), n \geq 0\}$\footnote{We normalize the Laplacian on $\mathbb{C}P^{N-1}$ by a factor $1/4$.}. Hence, the corresponding heat kernel is expanded in any orthonormal (with respect to the volume measure $\textrm{vol}_{\mathbb{C}P^{N-1}}$) basis of homogenous of degree zero spherical harmonics $(Y_j)_{j \geq 1}$ as: 
\begin{equation*}
R_t(w,z) \triangleq \sum_{n = 0}^{\infty}e^{-n(n+N-1)t} \sum_{j=1}^{d(n,N)}Y_j(w)\overline{Y_j}(z), \quad w,z\in \mathbb{C}P^{N-1}.
\end{equation*}
Here $d(n,N)$ is the dimension of the eigenspace of $(n,n)$-complex spherical harmonics given by (Theorem 3.6 in \cite{Koo})
\begin{equation*}
d(n,N) = \frac{2n+N-1}{N-1} \left(\frac{(N-1)_n}{n!}\right)^2.
\end{equation*}
Besides, the reproducing kernel formula (Theorem 3.8 in \cite{Koo}\footnote{The additional factor $2\pi$ comes from the fact that $\textrm{vol}(S^{2N-1}) = 2\pi \textrm{vol}(\mathbb{C}P^{N-1})$.}) 
shows that the kernel $R_t$ is real and does not depend on the choice of the basis (this is the analogue of (22) in \cite{Kar-McG}): 
\begin{equation}\label{analogue}
R_t(w,z) = 2\pi \sum_{n = 0}^{\infty}e^{-n(n+N-1)t}\frac{d(n,N)}{\textrm{vol}(S^{2N-1})} \frac{P_n^{N-2, 0}(2|\langle w,z\rangle|^2 - 1)}{P_n^{N-2, 0}(1)},
\end{equation}
where (\cite{And-Ask-Roy}, p.295)
\begin{equation*}
\langle w,z\rangle \triangleq \sum_{i=1}^N w_i\overline{z_i}, \quad P_n^{N-2,0}(1) = \frac{(N-1)_n}{n!},
\end{equation*}
and 
\begin{equation*}
\textrm{vol}(S^{2N-1}) = \frac{2\pi^N}{(N-1)!}.
\end{equation*}
is the volume of $S^{2N-1}$. Note that
\begin{equation*}
R_t(w,z) = \frac{(N-2)!}{\pi^{N-1}}\sum_{n = 0}^{\infty}e^{-n(n+N-1)t}\frac{P_n^{N-2, 0}(1) P_n^{N-2, 0}(2|\langle w,z\rangle|^2 - 1)}{||P_n^{N-2,0}||_2^2}
\end{equation*}
where 
\begin{equation*}
||P_n^{N-2,0}||_2^2 = \frac{1}{2n+N-1} 
\end{equation*}
is the squared $L^2$-norm of $u \mapsto P_n^{N-2,0}(2u-1)$ with respect to $s_1(u)du$ (\cite{And-Ask-Roy}, p.99). Thus Gasper's Theorem entails the positivity of $R_t$ (\cite{Gas}). Now, we proceed to the derivation of the density of $|U_t^1|^2$ and start with the decompositions 
\begin{eqnarray*}
w &=& \cos \theta_1 e_1 + \sin \theta_1 \xi_1, \\ 
z &=& \cos \theta_2 e_1 + \sin\theta_2 \xi_2 
\end{eqnarray*}
where $e_1$ is the first vector of the canonical basis of $\mathbb{C}^N$, $\theta_1,\theta_2 \in (0,\pi/2), \phi_1, \phi_2 \in (0,2\pi), \xi_1, \xi_2 \in S^{2N-3}$. The volume measure of $\mathbb{C} P^{N-1}$ in turn splits as (see \cite{Koo}, eq.2.18)
\begin{equation*}
\textrm{vol}_{\mathbb{C}P^{N-1}}(dz) = \cos\theta_2 (\sin\theta_2)^{2N-3} d\theta_2 \textrm{vol}_{S^{2N-3}}(d\xi_2). 
\end{equation*}
and the next step is to integrate \eqref{analogue} over $\xi_2$. But $U(N-1)$ acts transitively on $S^{2N-3}$ therefore we can take $\xi_1 = e_2$ to be the second vector of the canonical basis. As such, we are left with the volume of $S^{2N-5}$ (if $N$ is large enough) and with the integration over the distribution of the first coordinate of $\xi_2$. If this coordinate is parametrized by $(r,\psi)$ then its distribution reads 
\begin{equation*}
r(1-r^2)^{N-3}{\bf 1}_{[0,1]}(r){\bf 1}_{[0,2\pi]}(\psi) dr \, d\psi. 
\end{equation*}
Consequently, the density of $|U_t^{1}|^2$ displayed in \eqref{JacDen} follows from the product formula (4.12) in \cite{Koo} together with the variables change $u = \cos^2\theta_2$ ($c = \cos^2\theta_1$).

\begin{rem}
The eigenvalue of a $(n,n)$-spherical harmonic equals the eigenvalue of a $O(2N)$-spherical harmonic of degree $2n$ in $S^{2N-1}$ viewed as a real Euclidian sphere. This coincidence is due to the fact that both polynomials are homogenous with the same total degree $2n$ and since the correspond eigenvalue comes from the action of the Euler operator 
\begin{equation*}
\sum_{i=1}^N z_i \partial_{z_i} + \sum_{i=1}^N \overline{z_i}\partial_{\overline{z_i}} = \sum_{i=1}x_i\partial_{x_i} + \sum_{i=1}^Ny_i\partial_{y_i}, 
\end{equation*}
where $z_i \in \mathbb{C}$ is identified with $(x_i, y_i) \in \mathbb{R}^2$. 
\end{rem}


\section{The distribution of $(|U_t^1|^2, |U_t^2|^2)$}
Up to an additional ingredient, the lines of the previous proof enable to derive the density of $U_t^{(2)}$. More precisely, we start with the decompositions  
\begin{eqnarray*}
w &=& \cos \theta_1 e_1 + \sin \theta_1\xi_1 \\
& = & \cos \theta_1 e_1 + \sin \theta_1\cos \beta_1 e^{i\phi_1} e_2 + \sin\theta_1\sin \beta_1 \eta_1, \\ 
z &=& \cos \theta_2 e_1 + \sin\theta_2 \xi_2\\
&= &  \cos \theta_2 e_1 + \sin\theta_2\cos \beta_2 e^{i\phi_2} e_2 + \sin \theta_2\sin \beta_1 \eta_2,
\end{eqnarray*}
where $\beta_1, \beta_2 \in (0,\pi/2), \phi_1, \phi_2 \in (0,2\pi), \eta_1, \eta_2 \in S^{2N-3}$ and $e_2$ is the second vector of the canonical basis of $\mathbb{C}^N$. We also split the volume measure on $\mathbb{C}P^{N-1}$  as 
\begin{equation*}
\textrm{vol}_{\mathbb{C}P^{N-1}}(dz) = \left(\cos \theta_2 \sin^{2N-3}\theta_2 \cos\beta_2 \sin^{2N-5}\beta_2 d\theta_2 d\beta_2 d\phi_2\right) \textrm{vol}_{S^{2N-5}}(d\eta_2).
\end{equation*}
Now comes the needed additional ingredient, which is the special instance $m=n, \phi_1 = \phi_2 = 0$ in the formula stated in the bottom of p.5 in \cite{Koo1}. In order to recall it, let 
\begin{equation*}
p_j^{a,b}(x) \triangleq \frac{P_j^{a,b}(x)}{P_j^{a,b}(1)}, \quad a,b > -1
\end{equation*}
be the $j$-th normalized Jacobi polynomial and define the complex-valued polynomial (\cite{Koo} eq.3.15)
\begin{equation*}
R_{j,q}^{\alpha}(y) \triangleq |y|^{|j-q|}e^{i(j-q)\arg(y)} p_{j \wedge q}^{\alpha, |j-q|}(2|y|^2-1), \quad y \in \mathbb{C}, \, \alpha > -1.
\end{equation*}
as well as (\cite{Koo1} p.6)
\begin{equation*}
c_{j,q}(n,N) \triangleq \frac{N-2}{N-2+q+j} \binom{n}{q}\binom{n}{j} \frac{(N+n-1)_q (N+n-1)_j}{(N-2+j)_q (N-2+q)_j}.
\end{equation*}
Then the $n$-th reproducing kernel on $\mathbb{C}P^{N-1}$ admits the following expansion
\begin{align*}
\frac{P_n^{N-2, 0}(2|\langle w,z\rangle|^2 - 1)}{P_n^{N-2, 0}(1)} &= \sum_{j,q=0}^n c_{j,q}(n,N)[\sin(\theta_1)\sin\theta_2]^{j+q}[\cos\theta_1 \cos\theta_2]^{|j-q|}p_{n-j \wedge n-q}^{N-2+j+q, |j-q|}(\cos 2\theta_1)
\\& p_{n-j \wedge n-q}^{N-2+j+q, |j-q|}(\cos 2\theta_1) R_{j,q}^{N-3}\left(\langle \xi_1, \xi_2\rangle\right).
\end{align*}
 Substituting in \eqref{analogue}, we see that the next step towards the joint distribution of $(U_t^{1}, U_t^{2})$ consists in integrating 
\begin{equation*}
R_{j,q}^{N-3}\left(\langle \xi_1, \xi_2\rangle\right) = R_{j,q}^{N-3}\left(\cos \beta_1\cos\beta_2 + \sin\beta_1\sin \beta_2 \langle \eta_1, \eta_2\rangle\right)
\end{equation*}
over $\eta_2 \in S^{2N-3}$. To this end, we can assume without loss of generality that $\eta_1 = e_3$ (the third vector of the canonical basis) and use formula (4.11) in \cite{Koo}. Altogether, we get
\begin{align*}
&\frac{2\pi^2\textrm{vol}(S^{2N-5})}{(N-3)\textrm{vol}(S^{2N-1})}\sum_{n=0}^{\infty} e^{-n(n+N-1)t} d(n,N)\sum_{j,q=0}^n c_{j,q}(n,N)[\sin(\theta_1)\sin\theta_2]^{j+q}[\cos\theta_1 \cos\theta_2]^{|j-q|}
\\& p_{n-j \wedge n-q}^{N-2+j+q, |j-q|}(\cos 2\theta_1) p_{n-j \wedge n-q}^{N-2+j+q, |j-q|}(\cos 2\theta_2)R_{j,q}^{N-3}(\cos \beta_1 e^{i\phi_1})R_{j,q}^{N-3}(\cos \beta_2 e^{i\phi_2})
\\& = \frac{2(N-2)}{\pi}\sum_{n=0}^{\infty} e^{-n(n+N-1)t} (2n+N-1)[P_n^{N-2, 0}(1)]^2\sum_{j,q=0}^n c_{j,q}(n,N)[\sin(\theta_1)\sin\theta_2]^{j+q}[\cos\theta_1 \cos\theta_2]^{|j-q|}
\\& p_{n-j \wedge n-q}^{N-2+j+q, |j-q|}(\cos 2\theta_1)p_{n-j \wedge n-q}^{N-2+j+q, |j-q|}(\cos 2\theta_2)p_{j\wedge q}^{N-3, |j-q|}(\cos 2\beta_1)p_{j\wedge q}^{N-3, |j-q|}(\cos 2\beta_2) e^{i(j-q)(\phi_1+\phi_2)}.
 \end{align*}
with respect to 
\begin{equation*}
\cos \theta_2 \sin^{2N-3}\theta_2 \cos\beta_2 \sin^{2N-5}\beta_2 d\theta_2 d\beta_2d\phi_2.
\end{equation*}
Integrating over $\phi_2 \in (0,2\pi)$, then the sum over $(j, q)$ reduces to a sum over $q=j$. Thus, the density of $(\theta_2,\beta_2)$ given $(\theta_1, \beta_1)$ reads  
\begin{align*}
& 4(N-2)\sum_{n=0}^{\infty} e^{-n(n+N-1)t} (2n+N-1)[P_n^{N-2,0}(1)]^2 \sum_{j=0}^n c_{j,j}(n,N)[\sin(\theta_1)\sin\theta_2]^{2j}p_{n-j}^{N-2+2j,0}(\cos 2\theta_1)
\\& p_{n-j}^{N-2+2j, 0}(\cos 2\theta_2)p_{j}^{N-3, 0}(\cos 2\beta_1)p_{j}^{N-3, 0}(\cos 2\beta_2)
\end{align*}
with respect to 
\begin{equation*}
\cos \theta_2 \sin^{2N-3}\theta_2 \cos\beta_2 \sin^{2N-5}\beta_2 d\theta_2 d\beta_2.
\end{equation*}
Performing the variables change 
\begin{equation*}
u = \cos\theta_2, v = \sin\theta_2\cos\beta_2, 
\end{equation*}
we deduce that the density of $(|U_t^{1}|, |U_t^{2}|)$ given $(|U_0^{1}|, |U_0^{2}|)$ is:
\begin{align*}
& 4(N-2)\sum_{n=0}^{\infty}  e^{-n(n+N-1)t} (2n+N-1)[P_n^{N-2,0}(1)]^2 \sum_{j=0}^n c_{j,j}(n,N)[(1-|U_0^1|^2)(1-u^2)]^j 
\\& p_{n-j}^{N-2+2j,0}(2|U_0^1|^2-1) p_{n-j}^{N-2+2j, 0}(2u^2-1) p_{j}^{N-3, 0}\left(\frac{2|U_0^2|^2}{1-|U_0^1|^2} - 1\right)p_{j}^{N-3, 0}\left(\frac{2v^2}{1-u^2} - 1\right)
\end{align*}
 with respect to
\begin{equation*}
uv(1-u^2-v^2)^{N-3}{\bf 1}_{\{u > 0, v > 0, u^2+v^2 < 1\}}du\,dv.
\end{equation*}
Finally, if $(|U_0^{1}|^2, |U_0^{2}|^2) \triangleq (c_1,c_2) \in \Sigma_2$ then the density of $(|U_t^{1}|^2, |U_t^{2}|^2)$ reads
\begin{align*}
& (N-2)\sum_{n=0}^{\infty}  e^{-n(n+N-1)t} (2n+N-1)[P_n^{N-2,0}(1)]^2 \sum_{j=0}^n c_{j,j}(n,N)[(1-c_1)(1-u_1)]^j 
\\& p_{n-j}^{N-2+2j,0}(2c_1-1) p_{n-j}^{N-2+2j, 0}(2u_1-1) p_{j}^{N-3, 0}\left(\frac{2c_2}{1-2c_1} - 1\right)p_{j}^{N-3, 0}\left(\frac{2u_2}{1-u_1} - 1\right)
\end{align*}
 with respect to $s_2(u_1,u_2)du_1du_2$. The last expression may be put in a compact form as follows. For $n \geq 0$, set 
\begin{equation*}
Q_{n-j,j}^{(N)}(u,v) = (1-u_1)^jP_{n-j}^{N-2+2j, 0}(2u_1-1)P_{j}^{N-3, 0}\left(\frac{2u_2}{1-u_1} - 1\right), \quad j= 0, 1, \dots, n.
\end{equation*}
These are Jacobi polynomials in the simplex $\Sigma_2$ and are orthogonal with respect to the Dirichlet distribution whose density is $s_2(u_1,u_2)$ (specialize Proposition 2.3.8 in \cite{Dun-Xu} to $\alpha = (n-j, j), \kappa = (1/2, 1/2, N-5/2)$\footnote{Beware of the different normalization of the Jacobi polynomials used in Proposition 2.3.8. The reader is also invited to consult the list errata of available on the webpage of Y. Xu.}). After some computations, the density of $U_t^{(2)}$ may be written as 
\begin{equation*}
\sum_{n=0}^{\infty}  e^{-n(n+N-1)t} \sum_{j=0}^n\frac{Q_{n-j,j}^{(N)}(c_1,c_2)Q_{n-j,j}^{(N)}(u,v)}{||Q_{n-j,j}^{(N)}||_2^2}
\end{equation*}
where 
\begin{equation*}
||Q_{n-j,j}^{(N)}||_2^2 = \frac{1}{(2n+N-1)(2j+N-2)} = \frac{[P_{n-j}^{N-2+2j,0}(1)P_j^{N-3,0}(1)]^2}{(N-2)(2n+N-1)[P_n^{N-2,0}(1)]^2c_{j,j}(n,N)}
\end{equation*}
is the squared $L^2$-norm of $Q_{n-j,j}^{(N)}$ with respect to $s_2(u_1,u_2)du_1du_2$.
\begin{rem}
For general $k \geq 3$, the density of $U_t^{(k)}$ may be derived in a similar way by decomposing the variable $z\in \mathbb{C}P^{N-1}$ and the spherical harmonics on $S^{2N-1}$ over the sphere $S^{2N-2k+1}$. From the point of view of representation theory, this is equivalent to the decompositions of the representation of $\mathcal{U}(N)$ in the space of $\mathcal{U}(N)$-spherical harmonics under the action of the subgroup $\mathcal{U}(N-k+1)$. 
\end{rem}

\section{The distribution of $|U_t^1|^2$: another proof} 
In this section, we shall solve the system \eqref{S1} and prove \eqref{sequence}. To this end, we rewrite \eqref{S1} as,
\begin{equation}\label{S2}
\sum_{n=0}^p a_n \frac{1}{n!(p-n)!} \frac{1}{(N+2n)_{p-n}} = \frac{c^p}{(p!)^2}
\end{equation}
multiply both sides of \eqref{S2} by $(-1)^p(x/2)^{2p+N-1}$ for $x$ lying in some neighborhood of zero then sum over $p \geq 0$. Interchanging the order of summation, the system \eqref{S1} is equivalent to
\begin{equation*}
\sum_{n \geq 0}\frac{a_n}{n!}\Gamma(N+2n)(-1)^n J_{2n+N-1}(x) = J_0(\sqrt{c}x)\left(\frac{x}{2}\right)^{N-1}
\end{equation*}
where $J_{\alpha}$ is the Bessel function of index $\alpha \in \mathbb{R}$ defined by (\cite{Wat}): 
\begin{equation*}
J_{\alpha}(x) \triangleq \sum_{p \geq 0} \frac{(-1)^p}{p! \, \Gamma(p+\alpha+1)} \left(\frac{x}{2}\right)^{2p+\alpha}. 
\end{equation*}
Note in passing that the estimate (\cite{Wat})
\begin{equation*}
| \Gamma(N+2n)J_{2n+N-1}(x)| \leq \left(\frac{|x|}{2}\right)^{2n+N-1}
\end{equation*}
shows that \eqref{S2} converges provided 
\begin{equation}\label{S3}
\sum_{n \geq 0}\frac{|a_n|}{n!}\left(\frac{|x|}{2}\right)^{2n}
\end{equation}
does. Now recall the Neumann series (\cite{Wat}, p.138)
\begin{equation*}
\left(\frac{x}{2}\right)^{\nu} = \sum_{n \geq 0} \frac{(\nu+2n)\Gamma(\nu+n)}{n!} J_{\nu+2n}(x), \quad \nu \in \mathbb{N}. 
\end{equation*}
Specializing it to $\nu = N-2$, we get
\begin{align*}
(x/2)^{N-1}J_0(\sqrt{c}x) &= \sum_{p \geq 0} \frac{(-1)^pc^p}{(p!)^2}\sum_{n \geq 0} \frac{(2p+2n + N-1)\Gamma(2p+n + N-1)}{n!} J_{2p+2n+N-1}(x) 
\\& = \sum_{p \geq 0} \frac{(-1)^pc^p}{(p!)^2}\sum_{n \geq p} \frac{(2n + N-1)\Gamma(p+n+N-1)}{(n-p)!} J_{2n+N-1}(x) 
\\& = \sum_{n \geq 0} (2n+N-1)\sum_{p=0}^n \frac{(-1)^pc^p}{(p!)^2}  \frac{\Gamma(p+n+N-1)}{(n-p)!} J_{2n+N-1}(x). 
\end{align*}
Substituting in \eqref{S2}, then the uniqueness of the solution of $\eqref{S1}$ yields
\begin{equation*}
\frac{a_n}{n!}\Gamma(N+2n)(-1)^n =  (2n+N-1)\sum_{p=0}^n \frac{c^p}{(p!)^2}  \frac{\Gamma(p+n+N-1)}{(n-p)!}
\end{equation*}
or equivalently 
\begin{align*}
a_n & = \frac{(-1)^n n!}{\Gamma(N+2n-1)}\sum_{p=0}^n \frac{(-1)^pc^p}{(p!)^2}  \frac{\Gamma(p+n+N-1)}{(n-p)!} 
\\& = \frac{(-1)^n}{(N+n-1)_n}\sum_{p=0}^n \frac{(-1)^pn!}{(n-p)!}  \frac{(n+N-1)_p}{(p!)^2} c^p
\\& = \frac{(-1)^n}{(N+n-1)_n} {}_2F_1(-n, n+N-1, 1, c)
\end{align*}
where ${}_2F_1$ is the Gauss hypergeometric function (\cite{And-Ask-Roy}, p.). But from the very definition of Jacobi polynomials (\cite{And-Ask-Roy}, p.99)
\begin{equation*}
P_n^{\alpha, \beta}(x) = \frac{(\alpha+1)_n}{n!}{}_2F_1(-n, n+\alpha + \beta +1, \alpha+1, (1-x)/2), \,\, \alpha, \beta > -1,
\end{equation*}
and the relation $P_n^{\alpha,\beta} (x) = (-1)^n P_n^{\beta,\alpha}(-u)$ (\cite{And-Ask-Roy}, p.305), we obtain \eqref{sequence} as required. 
\begin{rem}
The estimate 
\begin{equation*}
|P_n^{0,N-2}(1-2c)| \leq P_n^{0,N-2}(1)
\end{equation*}
which follows for instance from the integral representation of the Gauss hypergeometric function shows that the series \eqref{S3} indeed converges absolutely everywhere. 
\end{rem}

With \eqref{sequence} in hands, we can invert the Laplace transform (\cite{Nec-Pel}, eq.4.11): 
\begin{equation*}
\varphi_{t/N}(c,\lambda) \triangleq \int_0^1 e^{\lambda u} f_{t/N}(c,u) du = \sum_{n = 0}^{\infty}a_ne^{-\Lambda_n t} \lambda^n {}_1F_1(n+1, N+2n, \lambda),\, \lambda \in \mathbb{R},
\end{equation*}
where ${}_1F_1$ is the confluent hypergeometric function (\cite{And-Ask-Roy}, \cite{Wat}) and $\Lambda_n = n(n+N-1)/N$. To proceed, recall the integral representation (\cite{And-Ask-Roy})
\begin{equation*}
{}_1F_1(a,b,\lambda) = \frac{\Gamma(b)}{\Gamma(b-a)\Gamma(a)} \int_0^1 e^{\lambda u}u^{a-1}(1-u)^{b-a-1}du, \quad b > a > 0. 
\end{equation*}
It follows that 
\begin{align*}
\lambda^n {}_1F_1(n+1, N+2n, \lambda) &= \frac{\Gamma(N+2n)}{\Gamma(N+n-1)n!} \int_0^1 \lambda^n e^{\lambda u}u^{n}(1-u)^{N+n-2}du
\\& = \frac{(-1)^n\Gamma(N+2n)}{\Gamma(N+n-1)n!}\int_0^1 e^{\lambda u}\left(\frac{d}{du}\right)^n [u^{n}(1-u)^{N+n-2}]du
\end{align*}
after $n$ integration by parts. But Rodriguez formula (\cite{And-Ask-Roy}, p.99)
\begin{equation*}
(1-x)^{\alpha}(1+x)^{\beta}P_n^{\alpha,\beta}(x) = \frac{(-1)^n}{2^n n!}\left(\frac{d}{du}\right)^n [(1-x)^{n+\alpha}(1+x)^{n+\beta}], \quad x \in (-1,1)
\end{equation*}
together with the variable change $x = 1-2u$ yields 
\begin{equation*}
\left(\frac{d}{du}\right)^n [u^{n}(1-u)^{N+n-2}] = n!(1-u)^{N-2}P_n^{0, N-2}(1-2u).
\end{equation*}
As a result 
\begin{align*}
a_n\lambda^n {}_1F_1(n+1, N+2n, \lambda) = P_n^{0,N-2}(1-2c) \int_0^1 e^{\lambda u} P_n^{0,N-2}(1-2u) \, (1-u)^{N-2} du
\end{align*}
and Tonelli-Fubini Theorem yields \eqref{JacDen} at time $t/N$. 

\section{From the Laplace transform to the generalized Jacobi operator} 
Another way to come from $\varphi_t(c,\lambda)$ to $f_t(c,\lambda)$ is as follows. For sake of simplicity, we shall drop the dependence on the parameter $c$. So, recall from \cite{Nec-Pel} Proposition 4.2 that $\varphi$ satisfies
\begin{equation*}
\partial_t\varphi = \lambda\varphi + \left(\lambda^2 - N\lambda\right)\partial_{\lambda}\varphi - \lambda^2\partial_{\lambda}^2\varphi
\end{equation*}
with the initial conditions $\varphi_0(c,\lambda) = e^{\lambda c}, \varphi_t(c,0) = 1$.
Assume  the density $f$ is unknown and is smooth in both variables $(t,u)$, then integration by parts yield
\begin{align*}
\int_0^1 e^{\lambda u}\partial_t f_t(u) du & = [e^{\lambda u}(1-Nu)f_t(u)]_0^1 + \int_0^1 e^{\lambda u}\partial_u[(Nu-1)f_t(u)]du 
\\& + [\lambda e^{\lambda u}u(1-u)f_t(u)]_0^1  - [e^{\lambda u} \partial_u(u(1-u)f_t(u))]_0^1 
\\& + \int_0^1 e^{\lambda u} \partial_u^2[u(1-u)f_t(u)] du
\\& = e^{\lambda}(2-N)f_t(1)  + \int_0^1 e^{\lambda u} \mathscr{L}(f_t)(u) du 
\end{align*}
where 
\begin{equation*}
\mathscr{L}_u \triangleq u(1-u)\partial_u^2 + [1+ (N-4)u]\partial_u + (N-2).
\end{equation*}
If $N = 2$ then 
\begin{equation*}
\mathscr{L} =  u(1-u)\partial_u^2 + [1- 2u]\partial_u
\end{equation*}
is nothing else but \eqref{JacOp} with $N=2$. Otherwise, write $f_t(u) \triangleq g_t(u) s_1(u)$ for a smooth function $g$ and note that $\mathscr{L}(s_1) = 0$. As a result,   
\begin{equation*}
\mathscr{L}(f_t)(u) = s_1(u)\left\{u(1-u)\partial_u^2 + [1- N u]\partial_u\right\}(g_t)(u) 
\end{equation*}
where the RHS is the operator displayed in \eqref{JacOp}. Since $f_t(1) = 0$ when $N \geq 3$ then we always have
\begin{equation*}
\int_0^1 e^{\lambda u}\partial_t g_t(u) s_1(u)du = \int_0^1 e^{\lambda u} \left\{u(1-u)\partial_u^2 + [1- N u]\partial_u\right\}(g_t)(u) s_1(u)du.
\end{equation*} 
But the set of monomials $(u^n)_{n \geq 0}$ is total in $L^2([0,1], (1-u)^{N-2}du)$ (\cite{Dun-Xu}, Theorem 3.17) then $g$ solves the heat equation
\begin{equation*}
\partial_t g_t(u) =  \left\{u(1-u)\partial_u^2 + [1- N u]\partial_u\right\}(g_t)(u). 
\end{equation*}
More generally, the Laplace transform of the density of $U_t^{(k)}, 1 \leq k \leq N-1$ satisfies the linear pde
\begin{eqnarray}\label{LinPde}
\partial_t\varphi&=&\sum_{j=1}^k \lambda_j \varphi+\sum_{j=1}^k(\lambda_j^2-N\lambda_j)\partial_j\varphi - \sum_{j,i=1}^k\lambda_i \lambda_j\partial_{ij}\varphi.
\end{eqnarray}
Hoping there will be no confusion, set again
$$
\varphi_t(c,\lambda) \triangleq \int_{\Sigma} e^{\langle \lambda,u\rangle }f_t(c,u)du
$$
where $du$ is the Lebesgue measure in the simplex $\Sigma$. Then, Integration by parts
\begin{eqnarray*}
\sum_{i=1}^k \lambda_i \varphi_t(\lambda)& \rightarrow &-\int_{\Sigma} e^{\langle \lambda,u\rangle}\left(\sum_{i=1}^k\partial_if_t(u)\right)du\\
(\lambda_i^2-N\lambda_i)\partial_i\varphi_t(\lambda)&\rightarrow &\int_{\Sigma} e^{\langle \lambda,u\rangle }[\partial_{ii}+N\partial_i](u_if_t(u))du, \quad i \in \{1,\dots, k\}, \\
\lambda_i^2\partial_{ii}\varphi_t(\lambda)&\rightarrow &\int_{\Sigma} e^{\langle \lambda,u\rangle }\partial_{ii}(u_i^2f_t(u))du, \quad i \in \{1,\dots, k\}, \\
\lambda_i\lambda_j\partial_{ij}\varphi_t(\lambda)&\rightarrow &\int_{\Sigma} e^{\langle \lambda,u\rangle }\partial_{ij}(u_i u_jf_t(u))du, \, \quad 1 \leq i \neq j \leq k,
\end{eqnarray*}
transform the pde \eqref{LinPde} into 
\begin{eqnarray*}
\int_{\Sigma} e^{\langle \lambda,u\rangle }\partial_t f_t(u) du = \textrm{boundary terms} + \int_{\Sigma} e^{\langle \lambda,u\rangle }\mathscr{L}(f_t)(u)du,
\end{eqnarray*}
where this time $\mathscr{L}$ denotes the operator
\begin{equation*}
 k(N- k-1)+\sum_{i=1}^k[1+[N- 4- 2(k-1)] u_i\big]\partial_i+ \sum_{i=1}^k(u_i-u_i^2)\partial_{ii} - \sum_{i \neq j}u_iu_j\partial_{ij}. 
\end{equation*}
If $k=N-1$ then $\mathscr{L}$ reduces to the operator displayed in \eqref{HE2}. Otherwise, set
\begin{equation*}
f_t(u) \triangleq g_t(u)s_k(u)
\end{equation*}
where this time $g$ is a smooth function in both variables $(t,u), u \in \Sigma_k$ and note that the relations $\partial_{i}s_k = \partial_{1}s_k, 1 \leq i \leq k$ together with the identity
\begin{equation*}
\sum_{i=1}^ku_i(1-u_i) - \sum_{i\neq j}u_iu_j = \sum_{i=1}^ku_i \left(1- \sum_{i=1}^ku_i \right) = \left(1- \sum_{i=1}^ku_i \right) - \left(1- \sum_{i=1}^ku_i \right)^2
\end{equation*}
imply that $\mathscr{L}_u(s_k) = 0$. Hence $\mathscr{L}(f_t)(u)$ gives rise to 
\begin{align*}
&\sum_{i=1}^k[1+[N- 4- 2(k-1)] u_i\big]\partial_i+ \sum_{i=1}^k(u_i-u_i^2)\partial_{ii} - \sum_{i \neq j}u_iu_j\partial_{ij} + 2\frac{\partial_1s}{s}(u) \left\{\sum_{i=1}^ku_i(1-u_i)\partial_i - \sum_{i \neq j}u_ju_i\partial_i\right\}
\\& = \sum_{i=1}^k[1+[N- 4- 2(k-1)] u_i\big]\partial_i+ \sum_{i=1}^k(u_i-u_i^2)\partial_{ii} - \sum_{i \neq j}u_iu_j\partial_{ij} + 2\frac{\partial_1s}{s}(u) \left\{\sum_{i=1}^ku_i\partial_i\left(1-u_i - \sum_{j \neq i}u_j\right)\right\}
\\&  = \sum_{i=1}^k[1 - N u_i]\partial_i+ \sum_{i=1}^k(u_i-u_i^2)\partial_{ii} - \sum_{i \neq j}u_iu_j\partial_{ij}
\end{align*}
acting on $g_t$. Consequently, if the boundary terms vanish then Theorem 3.17 in \cite{Dun-Xu} implies that $g_t$ solves the heat equation
\begin{equation*}
\left[\sum_{i=1}^k[1 - N u_i]\partial_i+ \sum_{i=1}^k(u_i-u_i^2)\partial_{ii} - \sum_{i \neq j}u_iu_j\partial_{ij}\right]g_t = \partial_tg_t.
\end{equation*}
We shall see below that this is the case provided that $1 \leq k \leq N-2$. 

\section{Analysis of the boundary terms}
Recall from the previous section that the integration by parts performed in the one-variable setting gave rise to the boundary term
\begin{equation*}
e^{\lambda u}\left[(1-Nu)f_t(u) -  \partial_u(u(1-u)f_t(u))\right]
\end{equation*}
which vanish at $u=1$ since $f_t(1) = 0$ when $N \geq 3$ (note that there is no such condition when $N=2$). For higher values $k \geq 2$, the situation is similar provided that $1 \leq k \leq N-2$ and is different when $k = N-1$ due to the interactions between $u_i$ and $u_j$ for $i \neq j$. Indeed, the boundary terms are given by  
\begin{align*}
& \sum_{i=1}^k \int [e^{\lambda_i u_i}(1-Nu_i)f_t(u)]_0^{1-\sum_{j\neq i}u_j}\left(\prod_{j\neq i} e^{\lambda_ju_j}{\bf 1}_{[0,1]}(u_j)du_j\right) + 
\\& \sum_{i=1}^k \int [\lambda_ie^{\lambda_i u_i}u_i(1-u_i)f_t(u)]_0^{1-\sum_{j\neq i}u_j}\left(\prod_{j\neq i} e^{\lambda_ju_j}{\bf 1}_{[0,1]}(u_j)du_j\right) - 
\\& \sum_{i=1}^k \int [e^{\lambda_i u_i}\partial_i\left\{u_i(1-u_i)f_t(u)\right\}]_0^{1-\sum_{j\neq i}u_j}\left(\prod_{j\neq i} e^{\lambda_ju_j}{\bf 1}_{[0,1]}(u_j)du_j\right) - 
\\& \sum_{1 \leq j\neq i \leq k}\int (\lambda_ju_j)[e^{\lambda_i u_i}u_if_t(u)]_0^{1-\sum_{m\neq i}u_m}\left(\prod_{m\neq i} e^{\lambda_mu_m}{\bf 1}_{[0,1]}(u_m)du_m\right) +
\\& \sum_{1 \leq j\neq i \leq k}\int [e^{\lambda_j u_j}u_j\partial_i\{u_if_t(u)\}]_0^{1-\sum_{m\neq j}u_m}\left(\prod_{m\neq j} e^{\lambda_mu_m}{\bf 1}_{[0,1]}(u_m)du_m\right).
\end{align*}
By Leibniz rule, the third and the last terms split into 
\begin{align*}
& \sum_{i=1}^k \int [e^{\lambda_i u_i}u_i(1-u_i)\partial_if_t(u)]_0^{1-\sum_{j\neq i}u_j}\left(\prod_{j\neq i} e^{\lambda_ju_j}{\bf 1}_{[0,1]}(u_j)du_j\right) + 
\\& \sum_{i=1}^k \int [e^{\lambda_i u_i}(1-2u_i)f_t(u)]_0^{1-\sum_{j\neq i}u_j}\left(\prod_{j\neq i} e^{\lambda_ju_j}{\bf 1}_{[0,1]}(u_j)du_j\right) 
\\& \end{align*}
and 
\begin{align*}
& \sum_{1 \leq j\neq i \leq k}\int u_i[e^{\lambda_j u_j}u_j\partial_if_t(u)]_0^{1-\sum_{m\neq j}u_m}\left(\prod_{m\neq j} e^{\lambda_mu_m}{\bf 1}_{[0,1]}(u_m)du_m\right) + 
\\&  (k-1)\sum_{j =1}^k\int [e^{\lambda_j u_j}u_jf_t(u)]_0^{1-\sum_{m\neq j}u_m}\left(\prod_{m\neq j} e^{\lambda_mu_m}{\bf 1}_{[0,1]}(u_m)du_m\right)
\end{align*}
respectively. Thus there are no boundary terms at $u_i=0, 1 \leq i \leq l$, while the remaining ones are given by 
\begin{align*}
& \sum_{i=1}^k \int [e^{\lambda_i u_i}(k+1-N)u_if_t(u)]^{u_i=1-\sum_{j\neq i}u_j}\left(\prod_{j\neq i} e^{\lambda_ju_j}{\bf 1}_{[0,1]}(u_j)du_j\right) +
\\& \sum_{i=1}^k \int \left[e^{\lambda_i u_i}u_i\left\{\lambda_i(1 -u_i) - \sum_{j \neq i}\lambda_ju_j\right\}f_t(u)\right]^{u_i=1-\sum_{j\neq i}u_j}\left(\prod_{j\neq i} e^{\lambda_ju_j}{\bf 1}_{[0,1]}(u_j)du_j\right) - 
\\& \sum_{i=1}^k\int \left[e^{\lambda_i u_i}u_i\left\{(1-u_i)\partial_if_t(u) - \sum_{j \neq i}u_j \partial_jf_t(u)\right\}\right]^{u_i=1-\sum_{j\neq i}u_j}\left(\prod_{j\neq i} e^{\lambda_ju_j}{\bf 1}_{[0,1]}(u_j)du_j\right).
\end{align*}
If $N \geq k+2$ and $f_t = g_t s_k$ vanishes on the hyperplane $\{u_1 + \dots + u_k = 1\}$ and the boundary terms reduce to 
\begin{equation*}
\sum_{i=1}^k\int e^{\lambda_i \left(1-\sum_{j\neq i}u_j\right)}\left(1-\sum_{j\neq i}u_j\right) \left\{\sum_{j\neq i}u_j[\partial_if_t(u) - \partial_jf_t(u)]\right\}^{u_i = 1-\sum_{j \neq i}u_j}\left(\prod_{j\neq i} e^{\lambda_ju_j}{\bf 1}_{[0,1]}(u_j)du_j\right).
\end{equation*}
But since $\partial_is_k = \partial_js_k$ and since $s_k$ vanishes on $\{u_1 + \dots + u_k = 1\}$, then for any $1 \leq i \neq j \leq k$
\begin{equation*}
\partial_if_t(u) = \partial_jf_t(u), \quad u_1 + \dots + u_k = 1
 \end{equation*}
so that all boundary terms vanish. When $k=N-1$ the boundary terms read 
\begin{align*}
& \sum_{i=1}^k \int e^{\lambda_i \left(1-\sum_{j\neq i}u_j\right)}\left(1-\sum_{j\neq i}u_j\right)  \left\{\sum_{j \neq i}(\lambda_i - \lambda_j)u_j\right\}\left[f_t(u)\right]^{u_i=1-\sum_{j\neq i}u_j}\left(\prod_{j\neq i} e^{\lambda_ju_j}{\bf 1}_{[0,1]}(u_j)du_j\right) - 
\\& \sum_{i=1}^k\int e^{\lambda_i \left(1-\sum_{j\neq i}u_j\right)}\left(1-\sum_{j\neq i}u_j\right) \left\{\sum_{j\neq i}u_j[\partial_if_t(u) - \partial_jf_t(u)]\right\}^{u_i = 1-\sum_{j \neq i}u_j}\left(\prod_{j\neq i} e^{\lambda_ju_j}{\bf 1}_{[0,1]}(u_j)du_j\right).
\end{align*}

{\bf Acknowledgments}: we would like to thank D. Bakry, C.F. Dunkl, T. Hmidi and Y. Xu for stimulating discussions and for their help with appropriate references.

\end{document}